\documentclass[11pt]{article}
\usepackage{amsmath}
\usepackage{latexsym}

\begin{document}

\title{Singular Vectors From Singular Values}
\author{
Weiwei Xu\thanks{Nanjing University of Information Science \& Technology, China (Email: wwx19840904@sina.com).
Research
supported in part by Natural Science Foundation of China under Grant No. 11971243 and
Natural Science Foundation of Jiangsu Province under grant BK20181405.}
\and
Michael K. Ng\thanks{
Department of Mathematics, The University of Hong Kong, Pokfulam, Hong Kong (Email: mng@maths.hku.hk).
Research
supported in part by the HKRGC GRF 12306616, 12200317, 12300218 and 12300519, and
HKU Grant 104005583.}}
\date{}
\maketitle

\begin{abstract}
\noindent
In the recent paper \cite{1}, Denton et al. provided the eigenvector-eigenvalue identity for Hermitian matrices,
and a survey was also given for such identity in the literature. The main aim of this paper is to present
the identity related to singular vectors and singular values of a general matrix.
\end{abstract}
\pagestyle{plain} \renewcommand{\thefootnote}{\fnsymbol{footnote}}



\section{Introduction}

In the recent paper \cite{1}, Denton et al. provided the eigenvector-eigenvalue identity for Hermitian matrices,
and a survey was also given for such identity in the literature. The main aim of this paper is to present
the identity related to singular vectors and singular values of a general matrix. Indeed, this topic has been studied by Thompson \cite{2}.
In particular, the interlacing inequalities for singular values of submstrices are derived in \cite{2}.

Let $A$ be a $m\times n$ matrix with singular values $\sigma_{i}(A)$ and normed left singular vectors $u_{i}$ and right singular vectors $v_{l}$.
$A^H$ denotes the conjugate transpose of $A$.
The elements of each left singular vector are denoted  $u_{ij}$ and right singular vector $v_{ls}$. Let $\bar{A}_{j}$ be $(n-1)\times n$ matrix after deleting the j-th row of $A$ and $\sigma_{k}(\bar{A}_{j})$ be singular values of $\bar{A}_{j}$. Let $\hat{A}_{s}$ be $n\times (n-1)$ matrix after deleting the s-th column of $A$ and $\sigma_{t}(\hat{A}_{s})$ be singular values of $\hat{A}_{s}$.

\textbf{Lemma 1}\cite{1}
\textit{
Let $A$ be a $n\times n$ Hermitian matrix with distinct eigenvalues $\lambda_{i}(A)$
and normed eigenvectors $v_{i}.$ The elements of each eigenvector are denoted $v_{i,j}$.
Let $M_{j}$ be the $(n-1)\times (n-1)$ submatrix of $A$ that results from deleting the j-th column and the j-th row, with eigenvalues $\lambda_{k}(M_{j})$. Then
\[
|v_{ij}|^{2}=\frac{\prod\limits_{k=1}^{n-1}(\lambda_{i}(A)-\lambda_{k}(M_{j}))}{\prod\limits_{k=1,k\neq i}^{n}(\lambda_{i}(A)-\lambda_{k}(A))},\;1\leq i,j\leq  n.
\]
}

\textbf{Lemma 2}\cite{1} \textit{
Let $A$ be a $n\times n$ Hermitian matrix with eigenvalues $\lambda_{i}(A)$
and normed eigenvectors $v_{i}.$ The elements of each eigenvector are denoted $v_{i,j}$.
Let $M_{j}$ be the $(n-1)\times (n-1)$ submatrix of $A$ that results from deleting the j-th column and the j-th row, with eigenvalues $\lambda_{k}(M_{j})$. Then
\[
|v_{ij}|^{2}\prod\limits_{k=1,k\neq i}^{n}(\lambda_{i}(A)-\lambda_{k}(A))=\prod\limits_{k=1}^{n-1}(\lambda_{i}(A)-\lambda_{k}(M_{j})),\;1\leq i,j\leq  n.
\]
}

In this paper, we establish the following theorem.
We first present singular vectors from distinct singular values of $n \times n$ matrix $A$.

\textbf{Theorem 1} \textit{Let $A=U^{H}\Sigma V\in\mathbf{C}^{n\times n}$ with distinct singular values $\sigma_{i}(A)$
and $U=(u_{ij}),1\leq i,j\leq n$ be left singular vector matrix of $A$ and $V=(v_{ls}),1\leq l,s\leq n$ be right singular
vector matrix of $A$. Let $\bar{A}_{j}$ be $(n-1)\times n$ matrix after deleting the j-th row of $A$ and $\sigma_{k}(\bar{A}_{j})$ be singular values of $\bar{A}_{j}$. Let $\hat{A}_{s}$ be $n\times (n-1)$ matrix after deleting the s-th column of $A$ and $\sigma_{t}(\hat{A}_{s})$ be singular values of $\hat{A}_{s}$. Then
\[
|u_{ij}|^{2}=\frac{\prod\limits_{k=1}^{n-1}(\sigma_{i}^{2}(A)-\sigma_{k}^{2}(\bar{A}_{j}))}{\prod\limits_{k=1,k\neq i}^{n}(\sigma_{i}^{2}(A)-\sigma_{k}^{2}(A))},\;1\leq i,j\leq  n,
\]
\[
|v_{ls}|^{2}=\frac{\prod\limits_{t=1}^{n-1}(\sigma_{l}^{2}(A)-\sigma_{t}^{2}(\hat{A}_{s}))}{\prod\limits_{t=1,t\neq l}^{n}(\sigma_{l}^{2}(A)-\sigma_{t}^{2}(A))},\;1\leq l,s\leq n.
\]
}

Then we deduce singular vectors from singular values of $m\times n$ matrix $A$.

\textbf{Theorem 2} \textit{Let $A=U^{H}\Sigma V\in\mathbf{C}_{r}^{m\times n}$ with singular values $\sigma_{1}(A)\geq\sigma_{2}(A)\geq\cdots\geq
\sigma_{r}(A)>0,\sigma_{r+1}(A)=\cdots=\sigma_{\min\{m,n\}}(A)=\cdots=\sigma_{\max\{m,n\}}(A)=0$
and $U=(u_{ij}),1\leq i,j\leq m$ be left singular vector matrix of $A$ and $V=(v_{ls}),1\leq l,s\leq n$ be right singular
vector matrix of $A$. Let $\bar{A}_{j}$ be $(m-1)\times n$ matrix after deleting the j-th row of $A$ and $\sigma_{k}(\bar{A}_{j})$ be singular values of $\bar{A}_{j}$. Let $\hat{A}_{s}$ be $m\times (n-1)$ matrix after deleting the s-th column of $A$ and $\sigma_{t}(\hat{A}_{s})$ be singular values of $\hat{A}_{s}$. Then
\[
|u_{ij}|^{2}\prod\limits_{k=1,k\neq i}^{m}(\sigma_{i}^{2}(A)-\sigma_{k}^{2}(A))=\prod\limits_{k=1}^{m-1}(\sigma_{i}^{2}(A)-\sigma_{k}^{2}(\bar{A}_{j})),\;1\leq i,j\leq  m,
\]
\[
|v_{ls}|^{2}\prod\limits_{t=1,t\neq l}^{n}(\sigma_{l}^{2}(A)-\sigma_{t}^{2}(A))=\prod\limits_{t=1}^{n-1}(\sigma_{l}^{2}(A)-\sigma_{t}^{2}(\hat{A}_{s})),\;1\leq l,s\leq n.
\]}

We note that the identities in Theorems 1 and 2 are not given in the previous results.
The organization of this paper is given as follows. In Section 2, we present the proofs of the two main theorems.
The concluding remarks are given in Section 3.

\section{The Proof}

In this section, we give the proofs of Theorems 1 and 2.

\noindent
\textit{Proof of Theorem 1.} Let singular value decompositions of $A$ be $A=U^{H}\Sigma V$, where $U\in\mathbf{U}^{n},V\in\mathbf{U}^{n}$.
(i) Let $A=\left(
\begin{array}{c}
a_{1}\\
a_{2}\\
\vdots\\
a_{n}
\end{array}
\right)$ with $a_{1},\ldots,a_{n}\in\mathbf{C}^{1\times n}.$ Since
\[
AA^{H}=U^{H}\Sigma\Sigma^{H}U,
\]
by Lemma 1 we have
\begin{eqnarray}
|u_{ij}|^{2}\prod\limits_{k=1,k\neq i}^{n}(\lambda_{i}(AA^{H})-\lambda_{k}(AA^{H}))=\prod\limits_{k=1}^{n-1}(\lambda_{i}(AA^{H})-\lambda_{k}(S_{j})),\;1\leq i,j\leq  n,
\end{eqnarray}
where $S_{j}$ is the $(n-1)\times (n-1)$ submatrix of $AA^{H}$ that results from deleting the $j$th column and the $j$th row, with eigenvalues $\lambda_{k}(S_{j})$. It is easy to see that
\[
S_{j}=\left(
\begin{array}{cccccc}
a_{1}a_{1}^{H}&\cdots &a_{1}a_{j-1}^{H}&a_{1}a_{j+1}^{H}&\cdots &a_{1}a_{n}^{H}\\
a_{2}a_{1}^{H}&\cdots &a_{2}a_{j-1}^{H}&a_{2}a_{j+1}^{H}&\cdots &a_{2}a_{n}^{H}\\
\cdots&\cdots&\cdots&\cdots&\cdots&\cdots\\
a_{j-1}a_{1}^{H}&\cdots &a_{j-1}a_{j-1}^{H}&a_{j-1}a_{j+1}^{H}&\cdots &a_{j-1}a_{n}^{H}\\
a_{j+1}a_{1}^{H}&\cdots &a_{j+1}a_{j-1}^{H}&a_{j+1}a_{j+1}^{H}&\cdots &a_{j+1}a_{n}^{H}\\
\cdots&\cdots&\cdots&\cdots&\cdots&\cdots\\
a_{n}a_{1}^{H}&\cdots &a_{n}a_{j-1}^{H}&a_{n}a_{j+1}^{H}&\cdots &a_{n}a_{n}^{H}
\end{array}
\right).
\]
Let $\bar{A}_{j}$ be $(n-1)\times n$ matrix after deleting the j-th row of $A$ and $\sigma_{k}(\bar{A}_{j})$ be singular values of $\bar{A}_{j}$, then
\[
\bar{A}_{j}=\left(
\begin{array}{c}
a_{1}\\
\vdots\\
a_{j-1}\\
a_{j+1}\\
\vdots\\
a_{n}
\end{array}
\right)
\]
and
\begin{eqnarray}
\sigma_{k}^{2}(\bar{A}_{j})=\lambda_{k}^{2}(\bar{A}_{j}\bar{A}_{j}^{H}).
\end{eqnarray}
Since $A$ has distinct singular values, by (1) and (2) we have
\[
|u_{ij}|^{2}=\frac{\prod\limits_{k=1}^{n-1}(\sigma_{i}^{2}(A)-\sigma_{k}^{2}(\bar{A}_{j}))}{\prod\limits_{k=1,k\neq i}^{n}(\sigma_{i}^{2}(A)-\sigma_{k}^{2}(A))},\;1\leq i,j\leq  n.
\]

(ii) Let $A=(b_{1},\ldots,b_{n})$ with $b_{1},\ldots,b_{n}\in\mathbf{C}^{n}.$ Since
\[
A^{H}A=V^{H}\Sigma^{H}\Sigma V,
\]
by Lemma 1 we have
\begin{eqnarray}
|v_{ij}|^{2}\prod\limits_{k=1,k\neq i}^{n}(\lambda_{i}(A^{H}A)-\lambda_{k}(A^{H}A))=\prod\limits_{k=1}^{n-1}(\lambda_{i}(A^{H}A)-\lambda_{k}(H_{j})),\;1\leq i,j\leq  n,
\end{eqnarray}
where $H_{j}$ is the $(n-1)\times (n-1)$ submatrix of $A^{H}A$ that results from deleting the $j$th column and the $j$th row, with eigenvalues $\lambda_{k}(H_{j})$. It is easy to see that
\[
H_{j}=\left(
\begin{array}{cccccc}
b_{1}^{H}b_{1}&\cdots &b_{1}^{H}b_{j-1}&b_{1}^{H}b_{j+1}&\cdots &b_{1}^{H}b_{n}\\
b_{2}^{H}b_{1}&\cdots &b_{2}^{H}b_{j-1}&b_{2}^{H}b_{j+1}&\cdots &b_{2}^{H}b_{n}\\
\cdots&\cdots&\cdots&\cdots&\cdots&\cdots\\
b_{j-1}^{H}b_{1}&\cdots &b_{j-1}^{H}b_{j-1}&b_{j-1}^{H}b_{j+1}&\cdots &b_{j-1}^{H}b_{n}\\
b_{j+1}^{H}b_{1}&\cdots &b_{j+1}^{H}b_{j-1}&b_{j+1}^{H}b_{j+1}&\cdots &b_{j+1}^{H}b_{n}\\
\cdots&\cdots&\cdots&\cdots&\cdots&\cdots\\
b_{n}^{H}b_{1}&\cdots &b_{n}^{H}b_{j-1}&b_{n}^{H}b_{j+1}&\cdots &b_{n}^{H}b_{n}
\end{array}
\right).
\]
Let $\hat{A}_{j}$ be $n\times (n-1)$ matrix after deleting the j-th row of $A$ and $\sigma_{k}(\hat{A}_{j})$ be singular values of $\hat{A}_{j}$, then
\[
\hat{A}_{j}=(b_{1},\ldots,b_{j-1},b_{j+1},\ldots,b_{n})
\]
and
\begin{eqnarray}
\sigma_{k}^{2}(\hat{A}_{j})=\lambda_{k}^{2}(\hat{A}_{j}^{H}\hat{A}_{j}).
\end{eqnarray}
Since $A$ has distinct singular values, by (3) and (4) we have
\[
|v_{ls}|^{2}=\frac{\prod\limits_{t=1}^{n-1}(\sigma_{l}^{2}(A)-\sigma_{t}^{2}(\hat{A}_{s}))}{\prod\limits_{t=1,t\neq l}^{n}(\sigma_{l}^{2}(A)-\sigma_{t}^{2}(A))},\;1\leq l,s\leq n.
\]
This completes the proof.

By using the similar trick, we can show the results of Theorem 2.

\noindent
\textit{Proof of Theorem 2.} Let singular value decompositions of $A$ be $A=U^{H}\Sigma V$, where $U\in\mathbf{U}^{m},V\in\mathbf{U}^{n}$.
(i) Let $A=\left(
\begin{array}{c}
a_{1}\\
a_{2}\\
\vdots\\
a_{m}
\end{array}
\right)$ with $a_{1},\ldots,a_{m}\in\mathbf{C}^{1\times n}.$ Since
\[
AA^{H}=U^{H}\Sigma\Sigma^{H}U,
\]
by Lemma 1 we have
\begin{eqnarray}
|u_{ij}|^{2}\prod\limits_{k=1,k\neq i}^{m}(\lambda_{i}(AA^{H})-\lambda_{k}(AA^{H}))=\prod\limits_{k=1}^{m-1}(\lambda_{i}(AA^{H})-\lambda_{k}(\tilde{S}_{j})),\;1\leq i,j\leq  m,
\end{eqnarray}
where $\tilde{S}_{j}$ is the $(m-1)\times (m-1)$ submatrix of $AA^{H}$ that results from deleting the $j$th column and the $j$th row, with eigenvalues $\lambda_{k}(\tilde{S}_{j})$. It is easy to see that
\[
\tilde{S}_{j}=\left(
\begin{array}{cccccc}
a_{1}a_{1}^{H}&\cdots &a_{1}a_{j-1}^{H}&a_{1}a_{j+1}^{H}&\cdots &a_{1}a_{m}^{H}\\
a_{2}a_{1}^{H}&\cdots &a_{2}a_{j-1}^{H}&a_{2}a_{j+1}^{H}&\cdots &a_{2}a_{m}^{H}\\
\cdots&\cdots&\cdots&\cdots&\cdots&\cdots\\
a_{j-1}a_{1}^{H}&\cdots &a_{j-1}a_{j-1}^{H}&a_{j-1}a_{j+1}^{H}&\cdots &a_{j-1}a_{m}^{H}\\
a_{j+1}a_{1}^{H}&\cdots &a_{j+1}a_{j-1}^{H}&a_{j+1}a_{j+1}^{H}&\cdots &a_{j+1}a_{m}^{H}\\
\cdots&\cdots&\cdots&\cdots&\cdots&\cdots\\
a_{m}a_{1}^{H}&\cdots &a_{m}a_{j-1}^{H}&a_{m}a_{j+1}^{H}&\cdots &a_{m}a_{m}^{H}
\end{array}
\right).
\]
Let $\tilde{A}_{j}$ be $(m-1)\times n$ matrix after deleting the j-th row of $A$ and $\sigma_{k}(\tilde{A}_{j})$ be singular values of $\tilde{A}_{j}$, then
\[
\tilde{A}_{j}=\left(
\begin{array}{c}
a_{1}\\
\vdots\\
a_{j-1}\\
a_{j+1}\\
\vdots\\
a_{m}
\end{array}
\right)
\]
and
\begin{eqnarray}
\sigma_{k}^{2}(\tilde{A}_{j})=\lambda_{k}^{2}(\tilde{A}_{j}\tilde{A}_{j}^{H}).
\end{eqnarray}
By (5) and (6) we have
\[
|u_{ij}|^{2}\prod\limits_{k=1,k\neq i}^{m}(\sigma_{i}^{2}(A)-\sigma_{k}^{2}(A))=\prod\limits_{k=1}^{m-1}(\sigma_{i}^{2}(A)-\sigma_{k}^{2}(\bar{A}_{j})),\;1\leq i,j\leq  m.
\]

(ii) Let $A=(b_{1},\ldots,b_{n})$ with $b_{1},\ldots,b_{n}\in\mathbf{C}^{m }.$ Since
\[
A^{H}A=V^{H}\Sigma^{H}\Sigma V,
\]
by Lemma 1 we have
\begin{eqnarray}
|v_{ij}|^{2}\prod\limits_{k=1,k\neq i}^{n}(\lambda_{i}(A^{H}A)-\lambda_{k}(A^{H}A))=\prod\limits_{k=1}^{n-1}(\lambda_{i}(A^{H}A)-\lambda_{k}(H_{j})),\;1\leq i,j\leq  n,
\end{eqnarray}
where $H_{j}$ is the $(n-1)\times (n-1)$ submatrix of $A^{H}A$ that results from deleting the $j$th column and the $j$th row, with eigenvalues $\lambda_{k}(H_{j})$. It is easy to see that
\[
H_{j}=\left(
\begin{array}{cccccc}
b_{1}^{H}b_{1}&\cdots &b_{1}^{H}b_{j-1}&b_{1}^{H}b_{j+1}&\cdots &b_{1}^{H}b_{n}\\
b_{2}^{H}b_{1}&\cdots &b_{2}^{H}b_{j-1}&b_{2}^{H}b_{j+1}&\cdots &b_{2}^{H}b_{n}\\
\cdots&\cdots&\cdots&\cdots&\cdots&\cdots\\
b_{j-1}^{H}b_{1}&\cdots &b_{j-1}^{H}b_{j-1}&b_{j-1}^{H}b_{j+1}&\cdots &b_{j-1}^{H}b_{n}\\
b_{j+1}^{H}b_{1}&\cdots &b_{j+1}^{H}b_{j-1}&b_{j+1}^{H}b_{j+1}&\cdots &b_{j+1}^{H}b_{n}\\
\cdots&\cdots&\cdots&\cdots&\cdots&\cdots\\
b_{n}^{H}b_{1}&\cdots &b_{n}^{H}b_{j-1}&b_{n}^{H}b_{j+1}&\cdots &b_{n}^{H}b_{n}
\end{array}
\right).
\]
Let $\hat{A}_{j}$ be $n\times (n-1)$ matrix after deleting the j-th row of $A$ and $\sigma_{k}(\hat{A}_{j})$ be singular values of $\hat{A}_{j}$, then
\[
\hat{A}_{j}=(b_{1},\ldots,b_{j-1},b_{j+1},\ldots,b_{n})
\]
and
\begin{eqnarray}
\sigma_{k}^{2}(\hat{A}_{j})=\lambda_{k}^{2}(\hat{A}_{j}^{H}\hat{A}_{j}).
\end{eqnarray}
By (7) and (8) we have
\[
|v_{ls}|^{2}=\frac{\prod\limits_{t=1}^{n-1}(\sigma_{l}^{2}(A)-\sigma_{t}^{2}(\hat{A}_{s}))}{\prod\limits_{t=1,t\neq l}^{n}(\sigma_{l}^{2}(A)-\sigma_{t}^{2}(A))},\;1\leq l,s\leq n.
\]
This completes the proof.

\section{Concluding Remarks}

In \cite{1}, Denton et al. revisited the identity of eigenvalue-eigenvector for Hermitian matrices.
This identity has been discovered by many researchers in the literature.
In this paper, we derived similar identity for singular vectors and singular values of general matrices.
However, we find these results have not been studied previously. As a future research work, it may be useful
to apply this new identity formula to matrix perturbation problems and
the calculation of arbitrary singular values and singular vector instead of all singular values and singular vectors.

\end{document}